\input amssym.def
\input amssym
\magnification=1200
\parindent0pt
\hsize=16 true cm
\baselineskip=13  pt plus .2pt

\def\S{\Bbb S}
\def\R{\Bbb R}
\def\Z{\Bbb Z}

\def\A{\Bbb A}

\centerline {\bf On minimal actions of finite simple groups on homology spheres}

\centerline {\bf  and Euclidean spaces}

\bigskip \bigskip

\centerline {Bruno P. Zimmermann}

\medskip

\centerline {Universit\`a degli Studi di Trieste}
\centerline {Dipartimento di Matematica e Informatica}
\centerline {34100 Trieste, Italy}

\bigskip \bigskip

{\it Abstract}.  We consider the following problem: for which classes of finite
groups, and in particular finite simple groups, does the minimal dimension of a
faithful, smooth action on a homology sphere coincide with the minimal
dimension of a faithful, linear action on a sphere?  We prove that the two
minimal dimensions coincide for the linear fractional groups ${\rm PSL}(2,p)$ as
well as for various classes of alternating and symmetric groups. We prove
analogous results also for actions on Euclidean spaces.

\bigskip  \bigskip

{\bf 1. Introduction}

\medskip

We consider actions of finite groups on spheres and homology spheres; if not stated
otherwise, all actions will be assumed to be smooth (or locally linear) and faithful.

\medskip

In a program to classify the finite, nonsolvable groups which
admit an action on a sphere or homology sphere in dimension three or four, it is
shown in [16-18], [20, 21] that the only finite (nonabelian) simple group acting on
a homology 3-sphere is the alternating or linear fractional group  $\A_5  \cong
{\rm PSL}(2,5)$, and the only finite simple groups acting on a homology 4-sphere
are the groups $\A_5 \cong {\rm PSL}(2,5)$ and $\A_6  \cong {\rm PSL}(2,9)$; a
corresponding classification in dimension five is still open.  This naturally
motivated the following problems:

\bigskip

-  For a given finite simple group $G$,  determine the minimal dimension of a sphere or homology
sphere on which $G$ acts faithfully (does it coincide with the minimal dimension of a linear
action  on a sphere?).

\medskip

-  Show that there are only finitely
many finite simple groups which admit an action on a homology
sphere of a given dimension $d$.

\bigskip

Some partial results for linear fractional and alternating groups have been
obtained in [22]. Concerning other classes of groups, it is shown in [9] that
for finite $p$-groups the two minimal dimensions of an action on a
homology sphere and of a linear action on a sphere coincide, and by [8] also
for orientation-preserving actions of finite abelian groups.
However, there are finite solvable groups for which the two minimal dimensions
do not coincide; specifically, the Milnor groups $Q(8a,b,c)$ ([15]) do not admit
faithful, orientation-preserving, linear actions on $S^3$ (neither free nor
nonfree), but by [14] some of them admit a free action on a homology 3-sphere
(see also section 4).

\medskip

In the present paper, we determine the minimal dimensions for the linear
fractional groups ${\rm PSL}(2,p)$ as well as for certain alternating and
symmetric groups. In the following, a mod $p$ homology sphere is a closed
$n$-manifold with the mod
$p$ homology of the $n$-sphere (i.e., homology with coefficients in the integers
$\Z_p$ mod $p$). Our main result is the following (where $p$ denotes a prime number and $q$
a positive integer):

\bigskip

{\bf Theorem 1.} {\sl  For the following finite groups, the minimal
dimension of a faithful, smooth (orientation-preserving) action on a mod $p$
homology sphere coincides with the minimal dimension of a faithful, linear
(orientation-preserving) action on a sphere:

\medskip

i) a semidirect product (metacyclic group)  $\Z_p \rtimes \Z_q$, with an
effective action of $\Z_q$ on the normal subgroup $\Z_p$;

\medskip

ii) a linear fractional group ${\rm PSL}(2,p)$;

\medskip

iii) a symmetric group $\S_p$;

\medskip

iv) an alternating groups $\A_p$ such that $p \equiv 3$  mod 4.}

\bigskip

See section 3 for the various minimal dimensions. An analogous result holds for
faithful, smooth actions on Euclidean spaces:

\bigskip

{\bf Theorem 2.} {\sl  For the groups listed in Theorem 1,
the minimal dimension of a faithful, smooth (orientation-preserving) action on
a Euclidean space coincides with the minimal dimension of a faithful,
linear  (orientation-preserving) action.}

\bigskip

Note that this is clear for actions with a global fixed point, so more
generally one may ask: what is the minimal dimension of a fixed point-free
action of a given finite group on a Euclidean space? See [7] for a survey on
actions of finite groups on disks and Euclidean spaces.

\bigskip

{\bf 2. The Borel spectral sequence associated to a group action}

\bigskip

Let $G$ be a finite group acting on a space $X$. Let $EG$ denote a contractible
space on which $G$ acts freely, and $BG = EG/G$.  We consider the twisted
product $X_G = EG \times_G X = (EG \times X)/G$. The "Borel fibering"  $X_G \to
BG$, with fiber $X$, is induced by the projection $EG \times X \to EG$, and the
equivariant cohomology of the $G$-space $X$ is defined as $H^*(X_G;K)$.  Our
main tool will be  the Leray-Serre spectral sequence $E(X)$  associated to the
Borel fibration $X_G \to BG$,

$$E_2^{i,j} = E_2^{i,j}(X) = H^i(BG;H^j(X;K)) = H^i(G;H^j(X;K)) \;\;\;
\Rightarrow
 \;\;\;  H^{i+j}(X_G;K),$$ i.e. converging to the graded group associated to a
filtration of $H^*(X_G;K)$ ("Borel spectral sequence", see e.g. [2]); here $K$
denotes any abelian coefficient group or commutative ring.

\medskip

Now suppose that $G$ acts on a (open or closed) $n$-manifold $M$; we denote by
$\Sigma$ the singular set of the $G$-action (all points in $M$ with nontrivial
stabilizer). Crucial for the proofs of Theorems 1 and 2 is the following
Proposition 1, see [13], [10, 11] for a proof (see also [4, Proposition VII.10.1]
for a Tate cohomology version ).

\bigskip

{\bf Proposition 1.} {\sl In dimensions greater than $n$, inclusion induces
isomorphisms $$H^*(M_G;K) \cong H^*(\Sigma_G;K).$$}

We use this to prove the following:

\bigskip

{\bf Proposition 2.} {\sl For an odd prime $p$ and a positive integer $q$, let
$G = \Z_p \rtimes \Z_q$ be a semidirect product with an effective action of
$\Z_q$ on the normal subgroup $\Z_p$.  Suppose that $G$ admits a faithful action
on a manifold $M$ with the mod $p$ homology of the $n$-sphere, and that the
subgroup $\Z_p$ of $G$ acts freely. Then $n+1$ is a multiple of $2q$ if all
elements of $G$ act as the identity on $H^n(M;\Bbb Z_p)\cong \Z_p$ (the
"orientation-preserving case"), and an odd multiple of $q$ if some element of
$G$ acts as the minus identity (the "orientation-reversing case").}

\medskip

{\it Proof.}  We consider first the Borel spectral sequence $E(\Sigma)$
converging to the cohomology of $\Sigma_G$, with $E_2^{i,j} =
H^i(G;H^j(\Sigma;\Bbb Z_p)$. Let $\Sigma_q$ denote the singular set of the
subgroup $\Z_q$ of $G$. The singular set $\Sigma$ of $G$ is the disjoint union
of the singular sets of the $p$ conjugates of $\Z_q$ in $G$; these fixed point
sets are all disjoint  since the action of $\Bbb Z_q$ on $\Z_p$ is effective and
the action of  $\Bbb Z_p$ is free, by assumption.  Now the action of $G$ on the
cohomology $H^j(\Sigma;\Bbb Z_p)$ of the orbit  $\Z_p(\Sigma_q) = \Sigma$ is
induced (or co-induced) from the action of $\Z_q$ on $H^j(\Sigma_q;\Bbb Z_p)$,
and by Shapiro's Lemma ([4, Proposition III.6.2]),
$H^i(G;H^j(\Sigma;\Bbb Z_p)$ is isomorphic to $H^i(\Z_q;H^j(\Sigma_q;\Bbb Z_p)$
and hence trivial, for $i > 0$. So also $H^*(\Sigma_G;\Bbb Z_p)$ is trivial, in
positive dimensions. By Proposition 1, also $H^*(M_G;\Z_p) \cong
H^*(\Sigma_G;\Z_p)$ is trivial, in sufficiently large dimensions.

\bigskip

Next we analyze the spectral sequence $E(M)$ converging to the cohomology of
$M_G$. The $E_2$-terms $E_2^{i,j} = H^i(G;H^j(M;\Bbb Z_p))$ are concentrated in
the two rows $j=0$ and $j=n$ where they are equal to $\Z_p$, with a possibly
twisted action of $G$ on $H^n(M;\Bbb Z_p) \cong \Z_p$. In particular, the only
possibly nontrivial differentials of $E(M)$ are  $d_{n+1}^{i,j}:E_{n+1}^{i,n}
\to E_{n+1}^{i+n+1,0}$, of bidegree $(n+1,-n)$.

\medskip

By [4, Theorem III.10.3], for $i>0$
$$H^i(G;H^j(M;\Bbb Z_p)) \cong H^i(\Z_p;H^j(M;\Bbb Z_p))^{\Z_q}.$$

The cohomology ring $H^*(\Z_p;\Z_p)$ is the tensor product of a polynomial
algebra $\Z_p[t]$ on a 2-dimensional generator $t$ and an exterior algebra
$E(s)$ on a 1-dimensional generator $s$ (see [1, Corollary II.4.2]); also, $t$
is the image of $s$ under the mod $p$ Bockstein homomorphism.

\medskip

Suppose first that $\Bbb Z_q$ acts trivially on $H^n(M;\Z_p) \cong \Z_p$. Since
the action of $\Z_q$ on $\Z_p$ is effective, also the action of $\Z_q$ on
$H^1(\Z_p;\Z_p)$ and $H^2(\Z_p;\Z_p)$ is effective: denoting by $\sigma$ the
action of a generator of $\Z_q$ on the cohomology $H^i(\Z_p;\Z_p)$, one has
$\sigma(s) = ks$  and $\sigma(t) = kt$, for some integer $k$ representing an
element of order $q$ in $\Z_p$. Hence $\sigma(t^a) = k^at^a$ and the only powers
of $t$ fixed by $\sigma$ are those divisible by $q$ (in dimensions $i$ which are
even multiples of $q$); similarly, $\sigma(st^a) = k^{a+1}t^{a+1}$, so the only
elements in odd dimensions fixed by $\sigma$ are the products
$st^a$ such that $a+1$ is a multiple of $q$ (in dimensions $i$ such that $i+1$
is an even multiple of $q$).   Consequently, $H^i(G;H^j(M;\Z_p))$ is nontrivial
exactly for $j=0$ and $n$ and in dimensions $i$ such that either $i$ or $i+1$ is
a multiple of $2q$.

\medskip

Since $H^*(M_G;\Z_p) \cong H^*(\Sigma_G;\Z_p)$ is trivial in sufficiently large
dimensions, the differentials $d_{n+1}$ of the spectral sequence $(E_2^{i,j}(M)
= H^i(G;H^j(M;\Z_p)$ (concentrated in the rows $j=0$ and $n$) have to be
isomorphisms in large dimensions. This can happen only if $n+1$ is a multiple of
$2q$, which completes the proof of Proposition 2 in the orientation-preserving
case.

\bigskip

Now suppose that a generator of $\Bbb Z_q$ acts as minus identity on
$H^n(M;\Z_p) \cong \Z_p$, in particular $q$ is even. Denoting by $\tilde
\sigma$ the action of a generator of $\Z_q$ on the cohomology
$H^i(\Z_p;H^n(M;\Z_p))$, we now have that  $\tilde \sigma(t^a) = -\sigma(t^a)=
-k^at^a$, $\tilde \sigma(st^a) = -\sigma(st^a) =  -k^{a+1}t^{a+1}$. Now the
elements of $H^i(\Z_p;H^n(M;\Z_p)) = H^i(\Z_p;\Z_p)$ invariant under $\tilde
\sigma$ are the powers of $t$ by odd multiples of $q/2$ (in dimensions $i$ which
are odd multiples of $q$), and the products $st^a$ such that $a+1$ is an odd
multiple of $q/2$ (in dimensions $i$ such that $i+1$ is an odd multiple of
$q$). Hence $H^i(G;H^j(M;\Z_p))$ is nontrivial exactly in the following
situations: either $j=0$ and $i$ or $i+1$ is an even multiple of $q$ (since the
action of $G$ on $H^0(M;\Z_p)$ is trivial), or if $j=n$ and $i$ or $i+1$ is an
odd multiple of $q$ (with the twisted action of $G$ on $H^n(M;\Z_p)$). Since
again $H^*(M_G;\Z_p) \cong H^*(\Sigma_G;\Z_p)$ is trivial in large dimensions,
the differential $d_{n+1}$ has to be an isomorphism and $n+1$ an odd multiple of
$q$.  (Note that, in order to obtain just the lower bound $n \ge q-1$, one may
apply the orientation-preserving case to the subgroup $\Z_p \rtimes
\Z_{q/2}$ of $\Z_p \rtimes \Z_q$.)

\medskip

This completes the proof of Proposition 2.

\bigskip

{\bf 3. Proofs of Theorems 1 and 2}

\bigskip

{\it Proof of Theorem 1.}

\medskip

i)  It is easy to see that a group $G = \Z_p \rtimes \Z_q$ as in Theorem 1 i)
admits a faithful, linear action (a faithful, real, linear representation) on
$S^{2q-1} \subset \R^{2q}$ if $q$ is odd, on $S^q \subset \R^{q+1}$ if $q$ is
even and the action is orientation-preserving, and on $S^{q-1} \subset \R^q$ if
some element of $G$ reverses the orientation (see e.g. [4, Example 9.2.3,
p.155]). We will show that the dimensions $2q-1$, $q$ and $q-1$ give lower bounds
for the dimensions of actions of $G$ on mod $p$ homology spheres, in the
respective cases.

\medskip

Suppose that $G = \Z_p \rtimes \Z_q$ admits a faithful, smooth action on a mod
$p$ homology $n$-sphere $M$. If the subgroup $\Bbb Z_p$ of $G$ acts freely,
Proposition 2 implies that $n+1$ is a multiple of $2q$ if $G$ acts
orientation-preservingly, and an odd multiple of $q$ otherwise. In particular,
the minimal possibilities are $n = 2q-1$ and $n = q-1$, respectively, which
coincide with the minimal dimension of a faithful, linear action on a sphere.

\medskip

Suppose then that $\Z_p$ has nonempty fixed point set $F$. Since the action is assumed to be
smooth, $F$ is a smooth submanifold; moreover, by Smith fixed point theory, $F$ is a mod $p$
homology sphere of some dimension $d$, $0 \le d  < n$. Then, by Lefschetz-duality with
coefficients in
$\Z_p$, the complement $M - F$ is a
$G$-invariant manifold with the mod $p$ homology of a sphere of dimension
$n-d-1$. Proposition 2 implies now that
$n-d-1 \ge 2q-1$, $n \ge 2q$ if $q$ is odd, $n-d-1 \ge q-1$, $n \ge q$
if $q$ is even.

\medskip

This completes the proof of part i) of Theorem 1.

\bigskip

ii)  The representation theory of the linear fractional groups ${\rm PSL}(2,p)$
is well-known, in particular ${\rm PSL}(2,p)$ admits a faithful, linear action
(a faithful, real, linear representation) on $S^{p-2} \subset \R^{p-1}$ if $p
\equiv 3$ mod 4, and on $S^{(p-1)/2} \subset \R^{(p+1)/2}$ if $p \equiv 1$ mod 4
(see e.g. [12], section 5.2, in particular Exercise 5.10).

\medskip

The subgroup of ${\rm PSL}(2,p)$ represented by all upper triangular matrices is
a semidirect product $\Z_p \rtimes \Z_q$, with $q=(p-1)/2$, where $\Z_q$ is the
subgroup represented by all diagonal matrices and $\Z_p$ by all matrices with
both diagonal entries equal to one; also, $\Z_q$ acts effectively on the normal
subgroup $\Z_p$. By i),  the minimal dimension of a faithful, smooth,
orientation-preserving action of $\Z_p \rtimes \Z_q$ on a mod $p$ homology
sphere is $2q-1 = p-2$ if $q$ is odd, and to $q = (p-1)/2$ if $q$ is even. Since
${\rm PSL}(2,p)$ admits linear actions in these dimensions, this proves Theorem
1 for the groups of type ii).

\bigskip

iii) and iv) Again the representation theory of the symmetric and alternating
groups is well-known (see e.g. [12, section 5.1]), in particular the symmetric
group $\S_p$ and the alternating group $\A_p$ admit a faithful, linear action (a
faithful, real, linear representation) on $S^{p-2} \subset \R^{p-1}$,
orientation-reversing and for an arbitrary positive integer $p$ in the case of
the symmetric group, for integers $p > 5$ in the case of the alternating group.

\medskip

For an odd prime $p$, consider the semidirect product $G = \Z_p \rtimes
\Z_{p-1}$, with an effective action of $\Z_{p-1}$ on the normal subgroup
$\Z_p$. The action by left-multiplication of $G$ on the $p$ left cosets of the
subgroup  $\Z_{p-1}$ of $G$ realizes $G$ as a subgroup of the symmetric group
$\S_p$. By i), the minimal dimension of a faithful, smooth action of $G$ on a
mod $p$ homology sphere is  $p-2$; since $\S_p$ admits a faithful, linear action
on $S^{p-2}$ this  proves part iii) of Theorem 1.

\medskip

For groups of type iv), we consider the subgroup $\Z_p \rtimes \Z_{(p-1)/2}$ of
index two of $\Z_p \rtimes \Z_{p-1}$ which is realized as a subgroup of the
alternating group $\A_p$. The alternating group $\A_p$ admits a faithful action
on  $S^{p-2}$; if $(p-1)/2$ is odd, by part i) of Theorem 1 this coincides with
the minimal dimension of a faithful, smooth action of $\Z_p \rtimes
\Z_{(p-1)/2}$ on a mod $p$ homology sphere.

\medskip

This completes the proof of Theorem 1.

\bigskip

{\it Proof of Theorem 2.}

\medskip

As noted in the proof of part i) of Theorem 1, the group $G = \Z_p \rtimes
\Z_q$ admits a faithful, linear action on $\R^{2q}$ if $q$ is odd, on
$\R^{q+1}$ if $q$ is even and the action is orientation-preserving, and on
$\R^q$ if some element of $G$ reverses the orientation.

\medskip

Suppose that $G = \Z_p \rtimes \Z_q$ admits a faithful, smooth action on
Euclidean space $\R^n$. By Smith fixed point theory (see [2]), the fixed point
set  $F$  of the subgroup $\Z_p$ of $G$ is a $\Z_p$-acyclic manifold of
some dimension $d \ge 0$, in particular non-empty. Note that the action of $G$
extends to a continuous action of $G$ on the sphere $S^n$ and, again by Smith
fixed point theory, the fixed point set $\hat F$ of the subgroup $\Z_p$ on
$S^n$ has the mod $p$ homology of a sphere of dimension $d$. Now by
Lefschetz-duality with coefficients in $\Z_p$, applied to the pair $(S^n, \hat
F)$, the complement $\R^n - F = S^n - \hat F$ has the mod $p$ homology of a
sphere of dimension $n - d - 1$. Hence $G$ admits a faithful, smooth action on
the manifold $\R^n - F$, with the mod $p$  homology of a sphere of dimension
$n - d - 1$, such that the normal subgroup $\Z_p$  of $G$ acts freely.  Now
Proposition 2 implies that $n - d - 1  \ge 2q - 1$, $n \ge 2q$ if $q$ is odd,
and $n - d - 1  \ge q - 1$, $n \ge q$ if $q$ is even.

\medskip

Suppose that $q$ is even and that $G$ acts orientation-preservingly on $\R^n$.
If the fixed point set $F$ of $\Z_p$ has dimension $d = 0$ then $F$ is a single
point which hence is a global fixed point of $G$. Now $G$ acts
orientation-preservingly on the boundary of an invariant neighbourhood of the
fixed point which is a sphere of dimension $n-1$, and Proposition 2 implies that
$n \ge 2q$. On the other hand, if $d > 0$ then, as noted above, $n-d-1 \ge q-1$
hence $n \ge q+d \ge q+1$.

\medskip

Since these lower bounds for the dimension of a faithful, smooth action of $G =
\Z_p \rtimes \Z_q$ on a Euclidean space can be realized by a faithful, linear
action of the group, this completes the proof of Theorem 2 for the groups of
type i). For all other groups, similar as in the proof of Theorem 1, it is a
consequence of the case of $\Z_p \rtimes \Z_q$.

\bigskip

{\bf 4.  Remarks on continuous actions of some other finite groups}

\bigskip

We present some other finite  groups for which the two minimal dimensions
for spheres and homology spheres coincide, even for the case of arbitrary
continuous actions.

\bigskip

{\bf Proposition 3.} {\sl  For the following groups, the minimal dimension of a
faithful, continuous action on a homology sphere
coincides with the minimal dimension of a faithful, linear action on a sphere:

\medskip

the unitary and symplectic groups $\; {\rm PSU}(4,2) \cong {\rm PSp}(4,3) \;$  and
$\; {\rm PSp}(6,2)$;

\medskip

the Weyl or Coxeter groups  ${\rm E}_6$, ${\rm E}_7$ and ${\rm E}_8$ of the
corresponding exceptional Lie algebras.}

\bigskip

This is a consequence of the following well-known result from Smith
fixed point theory ([19], see also [3]; iii) follows from the fact that the
fixed point set of an orientation-preserving involution on a manifold is a mod
2 homology manifold of codimension least two, then one applies ii)).

\bigskip

{\bf Proposition 4.}

\smallskip

{\sl  i) For an odd prime number $p$, the minimal
dimension of a faithful, continuous action of the elementary abelian $p$-group $(\Bbb
Z_p)^k$ on a mod $p$ homology sphere is $2k-1$ (and $2k$ for an action on a mod
$p$ acyclic manifold).

\medskip

ii)  The minimal dimension of a faithful, continuous action of the elementary
abelian 2-group $(\Bbb Z_2)^k$ on a mod 2  homology sphere is $k-1$ (and $k$ for an action on
a mod 2 acyclic manifold).

\medskip

iii)  The minimal dimension of a faithful, continuous, orientation-preserving
action of the elementary abelian 2-group $(\Bbb Z_2)^k$ on a mod 2 homology sphere is $k$
(and $k+1$ for an action on a mod 2 acyclic manifold).}

\bigskip

{\it Proof of Proposition 3.}    We refer to [6] and its references for
information about the subgroup structure and the character tables of the finite
simple groups.  The unitary group ${\rm PSU}(4,2)$ has a maximal subgroup $(\Bbb
Z_3)^3$ and a faithful, linear action on $\S^5$, so  Proposition 4 i) implies
that the two minimal dimensions coincide. Also, ${\rm PSU}(4,2)$ is a subgroup
of index two in the Weyl group of type ${\rm E}_6$ which has also a linear
action on $\S^5$.

\medskip

The Weyl group of type ${\rm E}_8$ has a subgroup $(\Bbb Z_3)^4$  and a linear
action on $\S^7$  (see the closely related  orthogonal group ${\rm O}^+_8(2)$ in
[6]), so the result follows again from  Proposition 4 i).

\medskip

The symplectic group  ${\rm PSp}(6,2)$ is a subgroup of index 2 in the Weyl
group of type ${\rm E}_7$, and both act linearly on $\S^6$.  Since ${\rm PSp}(6,2)$ has a
subgroup $(\Bbb Z_2)^6$, Proposition 4 iii) applies.  Alternatively, the group  ${\rm
PSp}(6,2)$ has the alternating group $\A_8$ and the linear
fractional group  ${\rm PSL}(2,8)$ as subgroups, and it
is shown in [22, Corollary 2] and [17, Proposition 1] as an application of the
Borel formula that the minimal dimension of a continuous action of
$\A_8$ and ${\rm PSL}(2,8)$ on a homology sphere is six (the results in
[22] and [17] are formulated for smooth actions; since the Borel formula holds for
continuous actions, the two specific results cited remain true in this more
general setting).

\medskip

This finishes the proof of Proposition 4.

\bigskip

As noted in the introduction, for some of the Milnor groups $Q(8a,b,c)$ the
minimal dimension of an action on a homology sphere is strictly smaller
than the minimal dimension of a linear action on a sphere. We do not
know similar examples for smooth or continuous actions on Euclidean space.
Interesting examples of continuous  actions not conjugate to smooth actions can
be obtained as follows.

\medskip

We consider again the Milnor groups $Q(8a,b,c)$ ([15]); these groups have
periodic cohomology of period four but do not admit a faithful, free, linear
action on the 3-sphere.  For odd, coprime integers $a > b > c \ge 1$, the Milnor
group
$Q(8a,b,c)$  is a semidirect product $\Z_{abc} \ltimes Q_8$ of a normal cyclic
subgroup  $\Z_a \times \Z_b \times \Z_c  \cong \Z_{abc}$ by the quaternion group
$Q_8 = \{\pm 1,\pm i,\pm j,\pm k\}$ of order eight, where $i,j$ and $k$ act
trivially on $\Z_a, \Z_b$ and $\Z_c$, respectively, and in a dihedral way on the
other two.

\medskip

It has been shown by Milgram [14] that some of the Milnor groups $Q(8a,b,c)$,
for odd, coprime integers $a > b > c \ge 1$, admit a faithful, free action on a
homology 3-sphere; let $Q$ be a Milnor group which admits such an action on a
homology 3-sphere $M$.  By the double suspension theorem (see e.g. [5]), the
double suspension $M * S^1$ of $M$ (or join with the 1-sphere) is homeomorphic
to $S^5$. Letting $Q$ act trivially on $S^1$, the actions of $Q$ on $M$
and $S^1$ induce a faithful, continuous, orientation-preserving action of $Q$ on
$S^5$ with fixed point set $S^1$, and hence also on $\R^5$ (the
complement of a fixed point).

\bigskip \bigskip

\centerline {\bf References}

\bigskip

\item {1.}  A.Adem, R.J.Milgram, {\it Cohomology of finite groups.}
Grundlehren der math. Wissenschaften 309, Springer-Verlag 1994

\item {2.} G.Bredon, {\it Introduction to compact Transformation Groups.}
Academic Press, New York 1972

\item {3.} M.R.Bridson, K.Vogtmann, {\it  Actions of automorphism groups of
free groups on homology spheres and acyclic manifolds.}  arXiv:0803.2062

\item {4.} K.S.Brown, {\it Cohomology of Groups.}  Graduate Texts in
Mathematics 87, Springer-Verlag 1982

\item {5.}  J.W.Cannon, {\it  The recognition problem: what is a topological
manifold.} Bull. Amer. Math. Soc.    , 832-866  (1978)

\item {6.} J.H.Conway, R.T.Curtis, S.P.Norton, R.A.Parker, R.A.Wilson, {\it
Atlas of Finite Groups.} Oxford University Press 1985

\item {7.}  M.W.Davis, {\it  A survey of results in higher dimensions.}  In:
The Smith Conjecture (eds. J.W.Morgan, H.Bass), Academic Press,  227-240  (1984)

\item {8.}  R.M.Dotzel, {\it  Orientation preserving actions of finite abelian
groups on spheres.}  Proc. Amer. Math. Soc. 100, 159-163 (1987)

\item {9.}  R.M.Dotzel, G.C.Hamrick, {\it  $p$-group actions on homology
spheres.}  Invent. math. 62, 437-442  (1981)

\item {10} A.E.Edmonds, {\it Aspects of group actions on four-manifolds.}
Top. Appl. 31, 109-124 (1989)

\item {11.} A.E.Edmonds, {\it Homologically trivial group actions on
4-manifolds.} Electronic version available at  arXiv:math.GT/9809055

\item {12.} W.Fulton, J.Harris, {\it Representation Theory: A First Course.}
Graduate Texts in Mathematics 129,  Springer-Verlag 1991

\item {13.} M.P.McCooey, {\it Symmetry groups of 4-manifolds.} Topology 41,
835-851  (2002)

\item {14.} R.J.Milgram, {\it Evaluating the Swan finiteness obstruction for
finite groups.} Algebraic and Geometric Topology. Lecture Notes in Math. 1126
(Springer 1985), 127-158

\item {15.} J.Milnor, {\it Groups which act on $S^n$ without fixed points.}
Amer. J. Math. 79, 623-630  (1957)

\item {16.} M.Mecchia, B.Zimmermann, {\it On finite groups acting on
$\Bbb Z_2$-homology 3-spheres.}  Math. Z. 248,   675-693  (2004)

\item {17.} M.Mecchia, B.Zimmermann, {\it On finite simple groups acting on
integer and mod 2 homology 3-spheres.}  J. Algebra 298, 460-467  (2006)

\item {18.} M.Mecchia, B.Zimmermann, {\it On finite simple and nonsolvable
groups acting on homology 4-spheres.}  Top. Appl. 2006, 2933-2942

\item {19.} P.A.Smith, {\it Permutable periodic transformations.} Proc. Nat.
Acad. Sci. U.S.A. 30, 105 - 108 (1944)

\item {20.} B.Zimmermann, {\it On finite simple groups acting on homology
3-spheres.}  Top. Appl. 125, 199-202 (2002)

\item {21.} B.Zimmermann, {\it On the classification of finite groups acting on
homology 3-spheres.} Pacific J. Math. 217, 387-395 (2004)

\item {22.} B.Zimmermann, {\it On the minimal dimension of a homology sphere on
which a finite group acts.}  Math. Proc. Camb. Phil. Soc. 144, 397-401 (2008)

\bye